\documentstyle[12pt]{article}
\newtheorem{theorem}{Theorem}[section]
\newtheorem{corollary}[theorem]{Corollary}
\newtheorem{definition}{Definition}
\newtheorem{example}[theorem]{Example}

\newtheorem{lemma}[theorem]{Lemma}
\title{Between ${\cal A}$- and ${\cal B}$-sets\thanks{1991 Math.\
Subject Classification --- Primary: 54G05, 54H05; Secondary:
54C05, 54C08, 54G99.\protect\newline Key words and phrases ---
${\cal A}$-set, ${\cal B}$-set, ${\cal A}{\cal B}$-set,
submaximal, partition space, extremally disconnected.
\protect\newline Research supported partially by the Ella and
Georg Ehrnrooth Foundation at Merita Bank, Finland.}}
\author{Julian Dontchev\\University of Helsinki\\Department of
Mathematics\\00014 Helsinki 10\\Finland}
\date{}
\begin{document}
\baselineskip=20pt plus 1pt minus 1pt
\newcommand{\abs}{\mbox{${\cal A}{\cal B}$-set}}
\newcommand{\abss}{\mbox{${\cal A}{\cal B}$-sets}}
\newcommand{\abc}{\mbox{${\cal A}{\cal B}$-continuous}}
\newcommand{\abco}{\mbox{${\cal A}{\cal B}$-continuity}}
\newcommand{\fxy}{\mbox{$f \colon (X,\tau) \rightarrow
(Y,\sigma)$}}
\maketitle
\begin{abstract} 
The aim of this paper is to introduce the class of \abss\ as the
sets that are the intersection of an open and a semi-regular set.
Several classes of well-known topological spaces are
characterized via the new concept. A new decomposition of
continuity is provided.
\end{abstract}

\section{Introduction}\label{s1}

A subset $A$ of a topological space $(X,\tau)$ is called {\em
locally closed} if $A$ is open in its closure or equivalently if
$A = U \cap V$, where $U$ is open and $V$ is closed. Several
classes of sets in general topological spaces have the above
mentioned property. For example, all connected subsets of the
real line as well as all locally compact subsets of Hausdorff
spaces are locally closed. Moreover, a Tychonoff topological
space $X$ is locally closed in its Stone-\v{C}ech
compactification ${\beta}X$ if and only if $X$ is locally
compact. Spaces in which every subset is locally closed are known
as submaximal. Recently, locally closed sets were studied in
\cite{GR0,GRV1}.

In 1986 and in 1989, Tong \cite{Tong1,Tong2} introduced two new
classes of set, namely $\cal A$-sets and $\cal B$-sets and using
them obtained new decompositions of continuity. He defined a set
$A$ to be an {\em $\cal A$-set} \cite{Tong1} (resp.\ a {\em $\cal
B$-set} \cite{Tong2}) if $A = U \cap V$, where $U$ is open and
$V$ is regular closed (resp.\ semi-closed (= t-set)). Clearly
every $\cal A$-set is locally closed and every locally closed set
is a $\cal B$-sets. Several topological spaces can be
characterized via the concepts of $\cal A$- and $\cal B$-sets
\cite{JD1}.

The concepts of $\cal A$-sets, locally closed sets and $\cal
B$-sets play important role when continuous functions are
decomposed. If the reader is interested in different
decompositions of continuity he (she) can refer to
\cite{ADG1,DG1,GGR1,GR1,GR2,Prz1,Tong1,Tong2}. Several new
decomposition of continuous and related mappings were recently
obtained in \cite{SDP1}.

The aim of this paper is to introduce a class of sets very
closely related to the classes of $\cal A$- and $\cal B$-sets,
in fact properly placed between them. Under consideration are the
sets that can be represented as the intersection of an open and
a semi-regular set. A subset $A$ of a (topological) space
$(X,\tau)$ is called {\em semi-regular} \cite{DMN1} if it is both
semi-open and semi-closed. In \cite{DMN1}, Di Maio and Noiri
pointed out that a set $A$ is semi-regular if and only if there
exists a regular open set $U$ such that $U \subseteq A \subseteq
\overline{U}$. Cameron \cite{C1} called semi-regular set {\em
regular semi-open}. In this paper, the connection of \abss\ to
other classes of `generalized open' sets is investigated as well
as several characterizations of topological spaces via \abss\ are
given. The concept of \abco\ is also introduced. A new
decomposition of continuity and a decomposition of \abco\ is
produced at the end of the paper.

Recall that a function $f \colon (X,\tau) \rightarrow (Y,\sigma)$
is called \"{A}-continuous \cite{Prz1} if for every open set $V$
of $(Y,\sigma)$, the set $f^{-1} (V)$ belongs to \"{A}, where
\"{A} is a collection of subsets of $X$. Most of the definitions
of function used throughout this paper are consequences of the
definition of \"{A}-continuity. However, for unknown concepts the
reader may refer to \cite{JD1,GR1}.

\section{${\cal A}{\cal B}$-sets}\label{s2}

\begin{definition}\label{d1}
{\em A subset $A$ of a space $(X,\tau)$ is called an {\em \abs\
} if $A = U \cap V$, where $U$ is open and $V$ is semi-regular.
The collection of all \abss\ in $X$ will be denoted by ${\cal
A}{\cal B} (X)$.}
\end{definition}

Since regular closed sets are semi-regular and since semi-regular
sets are semi-closed, then the following implications are
obvious:

\begin{center}
$\cal A$-set $\Rightarrow$ \abs\ $\Rightarrow$ $\cal B$-set
\end{center}

None of them of course is reversible as the following examples
shows:

\begin{example}\label{e1a}
{\em Let $X = \{ a,b,c,d \}$ and let $\tau = \{ \emptyset, \{ a
\}, \{ b \}, \{ a,b \}, X \}$. Set $A = \{ b,c \}$. It is easily
observed that $A$ is an \abs\ but not an $\cal A$-set.}
\end{example}

\begin{example}\label{e1b}
{\em Let $X$ be the space from Example 3.1 from \cite{Tong1},
i.e. let $X = \{ a,b,c \}$ and let $\tau = \{ \emptyset, \{ a \},
X \}$. Set $A = \{ c \}$. It is easily observed that $A$ is a
$\cal B$-set but not an \abs.}
\end{example}

Clearly every open and every semi-regular set is an \abs. But the
${\cal A}{\cal B}$-subset of the real line {\bf R} (with the
usual topology) $A = ({\bf R} \setminus \{ 0 \}) \cap [-1,1]$ is
neither open nor semi-regular.

Moreover, since the intersection of an open set and a
semi-regular set is always semi-open, then the following
implication is clear:

\begin{center}
\abs\ $\Rightarrow$ Semi-open set
\end{center}

However, if one considers the space from Example~\ref{e1b} above,
it becomes clear that not all semi-open sets are \abss: The set
$\{ a,b \}$ is semi-open but not an \abs.

Next the relation between $\cal B$-sets and \abss\ is shown but
first consider the following, probably known lemma. Recall that
a set $A \subseteq (X,\tau)$ is called {\em $\beta$-open} (=
semi-preopen) if $A \subset \overline{{\rm Int} \overline{A}}$.
The semi-closure of a set $A \subseteq (X,\tau)$ is the
intersection of all semi-closed supersets of $A$.

\begin{lemma}\label{l00}
The semi-closure of every $\beta$-open set is semi-regular.
$\Box$
\end{lemma}

\begin{theorem}\label{t00}
For a subset $A$ of a space $X$ the following are equivalent:

{\rm (1)} $A$ is an \abs.

{\rm (2)} $A$ is semi-open and a $\cal B$-set.

{\rm (3)} $A$ is $\beta$-open and a $\cal B$-set.
\end{theorem}

{\em Proof.} (1) $\Rightarrow$ (2) and (2) $\Rightarrow$ (3) are
obvious.

(3) $\Rightarrow$ (1) Since $A$ is a $\cal B$-set, then in the
notion of Theorem 1 from \cite{Y1}, there exists an open sets $U$
such that $A = U \cap {\rm sCl} A$, where ${\rm sCl} A$ denotes
the semi-closure of $A$ in $X$. By Lemma~\ref{l00}, ${\rm sCl}
A$ is semi-regular, since by (3) $A$ is $\beta$-open. Thus $A$
is an \abs. $\Box$

Recall that a space $X$ is called {\em submaximal} if every dense
subset of $X$ is open. Let ${\beta}O (X)$ denote the collection
of all $\beta$-open subset of $X$.

\begin{corollary}
If $(X,\tau)$ is a submaximal space, then ${\cal A}{\cal B} (X)
= {\beta}O (X)$.
\end{corollary}

{\em Proof.} Since $X$ is submaximal, then by Theorem 3.1 from
\cite{JD1} every ($\beta$-open) subset of $X$ is a $\cal B$-set.
Thus by Theorem~\ref{t00}, every $\beta$-open subset of $X$ is
an \abs. On the other hand, every \abs\ is $\beta$-open. $\Box$

The class of locally closed sets is also properly placed between
the classes of $\cal A$- and $\cal B$-sets but the concepts of
\abss\ and locally closed sets are independent from each other:
If first, every locally closed set is an \abs, then it would be
semi-open as well. But locally closed, semi-open sets are $\cal
A$-sets \cite[Theorem 1]{GR1}; however not all locally closed
sets are $\cal A$-sets. Second, if every \abs\ would be locally
closed, then again it must be an $\cal A$-set but as shown above
not all \abss\ are $\cal A$-sets.

\begin{theorem}\label{t0}
For a subset $A$ of a space $X$ the following are equivalent:

{\rm (1)} $A$ is semi-regular.

{\rm (2)} $A$ is semi-closed and an \abs.

{\rm (3)} $A$ is $\beta$-closed and an \abs.
\end{theorem}

{\em Proof.} (1) $\Rightarrow$ (2) and (2) $\Rightarrow$ (3) are
obvious.

(3) $\Rightarrow$ (1) Since $A$ is $\beta$-closed and a $\cal
B$-set, then $A$ is semi-closed \cite[Theorem 2.3]{JD1}. On the
other hand $A$ is semi-open, since it is an \abs. Thus $A$ is
semi-regular, being both semi-open and semi-closed. $\Box$

Recall that a subset $A$ of a space $(X,\tau)$ is called {\em
interior-closed} (= ic-set) \cite{GR2} if ${\rm Int} A$ is closed
in $A$. If $A \subseteq {\rm Int} \overline{A}$, then $A$ is
called {\em locally dense} \cite{CM1} (= preopen).

\begin{theorem}\label{t0a}
For a subset $A$ of a space $X$ the following are equivalent:

{\rm (1)} $A$ is open.

{\rm (2)} $A$ is an \abs\ and $A$ is either locally dense or an
ic-set.
\end{theorem}

{\em Proof.} (1) $\Rightarrow$ (2) is obvious.

(2) $\Rightarrow$ (1) If $A$ is locally dense, then since $A$ is
also a $\cal B$-set, it follows from Proposition 9 in
\cite{Tong2} that $A$ is open. If $A$ is an ic-set, then in the
notion of Theorem 1 from \cite{GR2}, $A$ is again open, since $A$
is also semi-open. $\Box$

\section{Some peculiar spaces}\label{s3}

Recall that a space $X$ is called {\em extremally disconnected}
(= ED) if every open subset of $X$ has open closure or
equivalently if every regular closed set is open.

\begin{theorem}\label{t1}
For a space $(X,\tau)$ the following are equivalent:

{\rm (1)} $X$ is ED.

{\rm (2)} $\tau = {\cal A}{\cal B} (X)$.

{\rm (3)} Every \abs\ is open.
\end{theorem}

{\em Proof.} (1) $\Rightarrow$ (2) Let $A \in {\cal A}{\cal B}
(X)$. Clearly $A$ is semi-open. From Theorem 4.1 in \cite{J1} it
follows that $A$ is preopen, since $X$ is ED. Moreover $A$ is a
$\cal B$-set and since it is preopen, it follows from Proposition
9 in \cite{Tong2} that $A \in \tau$. Hence ${\cal A}{\cal B} (X)
\subseteq \tau$. On the other hand it is obvious that $\tau
\subseteq {\cal A}{\cal B} (X)$.

(2) $\Rightarrow$ (3) is obvious.

(3) $\Rightarrow$ (1) Let $A \subseteq X$ be regular closed. Thus
$A$ is an \abs. By (3) $A$ is open. So, $X$ is ED. $\Box$

\begin{theorem}\label{t2}
For a space $X$ the following are equivalent:

{\rm (1)} $X$ is submaximal.

{\rm (2)} Every locally dense set is an \abs.

{\rm (3)} Every dense set is an \abs.
\end{theorem}

{\em Proof.} (1) $\Rightarrow$ (2) Let $A \subseteq X$ be locally
dense (= preopen). By (1), $A$ is open, since in submaximal
spaces every locally dense set is open \cite{JRV1}. Hence $A$ is
an \abs.

(2) $\Rightarrow$ (3) every dense set is locally dense.

(3) $\Rightarrow$ (1) Let $A \subseteq X$ be dense. By (3), $A$
is an \abs. Hence $A$ is both preopen and a $\cal B$-set. From
Proposition 9 in \cite{Tong2} it follows that $A$ is open. Thus
$X$ is submaximal. $\Box$

Recall that a space $X$ is called a {\em partition space} if
every open subset of $X$ is closed.

\begin{theorem}\label{t3}
For a space $X$ the following are equivalent:

{\rm (1)} $X$ is a partition space.

{\rm (2)} Every \abs\ is clopen.

{\rm (3)} Every \abs\ is (pre)closed.
\end{theorem}

{\em Proof.} (1) $\Rightarrow$ (2) Let $A \subseteq X$ be an
\abs. By (1) and Theorem 3.2 from \cite{JD1}, $A$ is clopen,
since it is a $\cal B$-set.

(2) $\Rightarrow$ (3) every clopen set is preclosed.

(3) $\Rightarrow$ (1) Let $A \subseteq X$ be open. Then $A$ is
an \abs\ and by (3) it is preclosed. Since every preclosed
(semi-)open set is (regular) closed, then $X$ is a partition
space. $\Box$

\begin{theorem}\label{t4}
For a space $X$ the following are equivalent:

{\rm (1)} $X$ is indiscrete.

{\rm (2)} ${\cal A}{\cal B} (X) = \{ \emptyset, X \}$.
\end{theorem}

{\em Proof.} The theorem follows from Theorem 3.3 from
\cite{JD1}, since the class of \abss\ is (properly) placed
between
the classes of $\cal A$- and $\cal B$-sets. $\Box$

\begin{theorem}\label{t5}
For a space $X$ the following are equivalent:

{\rm (1)} $X$ is discrete.

{\rm (2)} Every subset of $X$ is an \abs. 

{\rm (3)} Every singleton is an \abs.
\end{theorem}

{\em Proof.} (1) $\Rightarrow$ (2) and (2) $\Rightarrow$ (3) are
obvious.

(3) $\Rightarrow$ (1) Let $x \in X$. By (3), $\{ x \}$ is an
\abs\ and hence semi-open. Then $\{ x \}$ must contain a non-void
open subset. Since the only possibility is $\{ x \}$ itself, then
each singleton is open or equivalently $X$ is discrete. $\Box$

Recall that a space $X$ is called {\em hyperconnected} if every
open subset of $X$ is dense in $X$.

\begin{theorem}\label{t6}
For a space $X$ the following are equivalent:

{\rm (1)} $X$ is hyperconnected.

{\rm (2)} Every \abs\ is dense.
\end{theorem}

{\em Proof.} (1) $\Rightarrow$ (2) Let $A \subseteq X$ be an
\abs. Then $A$ is semi-open and hence there exist an open subset
$U$ such that $U \subseteq A \subseteq \overline{U}$. By (1), $U$
is dense. Hence its superset $A$ is also dense.

(2) $\Rightarrow$ (1) Every open subset of $X$ is an \abs\ and
hence by (2) dense. $\Box$

Recall that a space $X$ is called {\em semi-connected} \cite{PR1}
if $X$ cannot be expressed as the disjoint union of two non-void
semi-open sets.

\begin{theorem}\label{t7}
For a space $X$ the following are equivalent:

{\rm (1)} $X$ is semi-connected.

{\rm (2)} $X$ is not the union of two disjoint non-void \abss.
\end{theorem}

{\em Proof.} (1) $\Rightarrow$ (2) If $X$ is the union of two
disjoint non-void \abss, then $X$ is not semi-connected, since
\abss\ are semi-open.

(2) $\Rightarrow$ (1) If $X$ is not semi-connected, then $X$ has
a non-trivial semi-open subset $A$ with semi-open complement.
Since both $A$ and $B = X \setminus A$ are semi-regular, then $A$
and $B$ are \abss. So $X$ is the union of two disjoint non-void
\abss, contradictory to (2). $\Box$

\section{\abc\ functions}\label{s4}

\begin{definition}\label{d2}
{\em A function \fxy\ is called {\em \abc\ } if the preimage of
every open subset of $Y$ is an \abs\ in $X$.}
\end{definition}

Recall that a function \fxy\ is called {\em strongly irresolute}
\cite{DCP1} if $f({\rm sCl} A) \subseteq f(A)$ for every subset
$A$ of $X$. It is easily observed that a function \fxy\ is
strongly irresolute if and only if the inverse image of every
subset of $Y$ is semi-regular in $X$.

The last four theorems are consequences of results from the
beginning of this paper, therefore their proofs are omitted.
Theorem~\ref{s41} gives the relations between \abc\ functions and
other forms of `generalized continuity'. Note that none of the
implications in Theorem~\ref{s41} is reversible.
Theorem~\ref{s42} gives a decomposition of \abco, while
Theorem~\ref{s42a} is an improvement of Theorem 4 (i) from
\cite{GR1} and it follows from Theorem 2.4 in \cite{JD1}.
Theorem~\ref{s43} gives a decomposition of continuity dual to
\abco.

\begin{theorem}\label{s41}
{\rm (i)} Every $\cal A$-continuous function is \abc,

{\rm (ii)} Every strongly irresolute function is \abc,

{\rm (iii)} Every \abc\ function is $\cal B$-continuous,

{\rm (iv)} Every \abc\ function is semi-continuous. $\Box$
\end{theorem}

\begin{theorem}\label{s42}
For a function \fxy, the following conditions are equivalent:

{\rm (1)} $f$ is \abc.

{\rm (2)} $f$ is semi-continuous and $\cal B$-continuous.

{\rm (3)} $f$ is $\beta$-continuous and $\cal B$-continuous.
$\Box$
\end{theorem}

\begin{theorem}\label{s42a}
For a function \fxy, the following conditions are equivalent:

{\rm (1)} $f$ is $\cal A$-continuous.

{\rm (2)} $f$ is $\beta$-continuous and LC-continuous. $\Box$
\end{theorem}

\begin{theorem}\label{s43}
For a function \fxy, the following conditions are equivalent:

{\rm (1)} $f$ is continuous.

{\rm (2)} $f$ is \abc\ and either precontinuous or ic-continuous.
$\Box$
\end{theorem}

\
\
e-mail: {\tt dontchev@cc.helsinki.fi}
\
\
\end{document}